\documentclass[12pt,reqno]{amsart} 
\usepackage{amsmath}
\usepackage{amssymb, latexsym, amsfonts, amscd, amsthm, mathrsfs, enumerate, esint}
\usepackage[usenames,dvipsnames]{color}
\usepackage[all]{xy}
\usepackage{graphicx}
\usepackage{epsfig}
\usepackage{hyperref}

\numberwithin{equation}{section}

\definecolor{purple}{rgb}{0.9,0,0.8}

\definecolor{gray}{rgb}{0.7,0.7,0.7}

\newcommand{\abbr}[1]{{\sc\lowercase{#1}}}

\newtheorem{thm}{Theorem}[section]

\newtheorem{ppn}[thm]{Proposition}

\newtheorem{defn}[thm]{Definition}
\theoremstyle{definition}

\newtheorem{remark}[thm]{Remark}
\newcommand{\beq}{\begin{equation}}
\newcommand{\eeq}{\end{equation}}

\newcommand{\ep}{\epsilon} 
\newcommand{\vep}{\varepsilon}

\newcommand{\bB}{\mathbb{B}}

\newcommand{\bG}{\mathbb{G}}

\newcommand{\bI}{\mathbb{I}}

\newcommand{\bN}{\mathbb{N}}

\newcommand{\bP}{\mathbb{P}}

\newcommand{\bR}{\mathbb{R}}

\newcommand{\V}{\mathbb{V}}

\newcommand{\bZ}{\mathbb{Z}}

\newcommand{\cE}{\mathcal{E}}
\newcommand{\cF}{\mathcal{F}}

\newcommand{\cL}{\mathcal{L}}

\newcommand{\cS}{\mathcal{S}}

\newcommand{\cZ}{\mathcal{Z}}

\newcommand{\sC}{\mathscr{C}}

\newcommand{\wh}{\widehat}

\newcommand{\grad}{\nabla}

\begin{document}

\title[Stability and instability of Gaussian heat kernel estimates]
{Stability and instability of Gaussian heat kernel estimates for random walks among 
time-dependent conductances}
\date{\today}

\author[R.\ Huang]{Ruojun Huang$^*$}
\author[T.\ Kumagai]{Takashi Kumagai$^\diamond$}
\address{$^*$Department of Statistics, Stanford University,  
Sequoia Hall, 390 Serra Mall, 
Stanford, CA 94305, USA. E-mail: {\tt hruojun@stanford.edu}}
\address{$^{\diamond}$Research Institute for Mathematical Sciences,
Kyoto University, Kyoto 606-8502, Japan. E-mail: {\tt kumagai@kurims.kyoto-u.ac.jp}}

\thanks{This research was supported in part 
by JSPS KAKENHI Grant Number 25247007 and by NSF grant DMS-1106627.}



\maketitle
\begin{abstract} We consider time-dependent random walks among time-dependent conductances. 
For discrete time random walks, we show that, unlike the time-independent case, two-sided Gaussian heat kernel estimates are not stable under perturbations. This is proved by giving an example 
of a ballistic and transient time-dependent random walk on $\bZ$ among uniformly elliptic time-dependent 
conductances. For continuous time random walks, we show the instability when the holding times are 
i.i.d. $\exp(1)$, and in contrast, we prove the stability when the holding times change by sites in such a way that the base measure is a uniform measure.
\end{abstract}

\begin{section}{Introduction}\label{statement}
The study of heat kernels of diffusions on manifolds and Markov chains on graphs has a very long and fruitful history. One of the motivations was to obtain a priori estimates such as the estimates
of the H\"older continuity for the solutions of heat equations. 
In the framework of the divergence operator 
${\mathcal L}=\sum_{i,j=1}^d\frac{\partial}{\partial x_i}
(a_{ij}(x)\frac{\partial}{\partial x_j})$ 
on ${\mathbb R}^d$ where $a_{ij}(\cdot)$ is measurable and 
symmetric, there are significant work by De Giorgi, Nash and Moser around late 50s to early 60s. 
For the divergence form satisfying a uniformly elliptic
condition,  Aronson \cite{Ar} proved the following two-sided Gaussian heat kernel estimates
for all $t>0, x,y\in{\mathbb R}^d$:
\begin{equation}
c_1t^{-d/2}\exp \Big(- \frac {c_2d(x,y)^{2}}{t}\Big)
\le p_t(x,y)\le c_3t^{-d/2}\exp\Big(- \frac {c_4d(x,y)^{2}}{t}\Big).
\label{eq:aron}\end{equation}
Later in the last century, the two-sided Gaussian estimates were obtained for many operators in many spaces and the heat kernel estimates were investigated from various aspects. One of the important directions is to establish the stability of the estimates, namely
to show that the estimates are preserved when the operator (or the corresponding Dirichlet form) is perturbed in a suitable way.  Consider the Laplace-Beltrami operator on a complete Riemannian manifold
with $d(\cdot,\cdot)$ and $\mu$ being the Riemannian metric and the Riemannian measure. 
Early in 90s, Grigor'yan \cite{Gr} and Saloff-Coste \cite{SC} independently proved that for the Laplace-Beltrami operator,    
a variant of \eqref{eq:aron} (i.e. changing $t^{-d/2}$ into $\mu(\bB(x,t^{1/2}))^{-1}$) 
 is equivalent to a volume doubling condition (VD)
plus Poincar\'e inequalities (PI(2)) via the equivalence to 
parabolic Harnack inequalities --see Theorem \ref{thm:PIVDPHI} for definitions of the terminologies. 
Since (VD) and (PI(2)) are 
stable under the perturbations, one can obtain the stability of the heat kernel estimates. 
The results were later extended to the framework of
Dirichlet forms on metric measure spaces by Sturm \cite{St1,St2} and graphs by Delmotte \cite{De}.  
We note that such a stability theory has been extended to the sub-Gaussian heat kernel estimates, 
also to locally irregular graphs such as the super-critical Bernoulli percolation cluster (\cite{BC});  
the theory has been very useful in the recent developments of the random walk among random conductances (see for example \cite{Bisk,Kum}). 


In this note, 
we are mainly interested in cases when the edge conductances of the graph are themselves changing in time, independently of the walk. 
We will consider the stability of the two-sided Gaussian heat kernel estimates in this setting. 
A naive guess is that the stability holds at least 
when the time-dependent conductance is bounded from above and below uniformly by positive constants. 
However, this naive guess is completely wrong for discrete time (time-dependent) random walks. 
Indeed, in Proposition \ref{ctr_ex}(i), we give an example of a ballistic and transient time-dependent random walk 
on $\bZ$ among uniformly elliptic time-dependent conductances. We also give a counter example in the setting  of continuous time (time-dependent) random walks, called constant speed random walks, when the holding times are i.i.d. $\exp(1)$ (Proposition \ref{ctr_ex}(ii)). Contrary to the above, when the holding times change by sites in such a way that the base measure is a uniform measure (called a variable speed random walk), we can prove the stability by proving the equivalence of the heat kernel estimates to 
(VD) and (PI(2)) under some regularity condition of the conductances --see Theorem \ref{thm:PIVDPHI}. We note that the stability 
of parabolic Harnack inequalities and estimates of the heat kernel were already established 
in the framework of time-dependent Dirichlet forms on metric measure spaces by Sturm \cite{St1, St2}, and in \cite{DD, GOS}
it was proved that for random walks on ${\mathbb Z}^d$ among uniformly elliptic time-dependent conductances, the two-sided Gaussian heat kernel estimates hold. 
Also in \cite{GP}, some criteria was given on the recurrence and transience of a set using the heat kernel estimates. These are results for variable speed random walks.
The purpose of this note is to demonstrate a fundamental difference between discrete time  
random walks (or constant speed random walks) and variable speed random walks even in the framework 
that the random walks are uniformly elliptic and uniformly lazy.   
In particular, we find that 
both the upper and lower Gaussian bounds can be violated in these situations (see Proposition 
\ref{ctr_ex} and \ref{ctr_ex_2}), 
and consequently they are unstable even though (VD)+(PI(2)) are still stable.
This contrasts to the above mentioned 
situation on graphs with time-independent conductances, where the discrete time and the two types of continuous time random walks share the same long-time properties at least in the uniformly elliptic setting. 

Let us mention some related works. \cite{ABGK,DHS} study recurrence versus transience of discrete time simple random walks on graphs with monotonically changing conductances. 
To be fair, we note 
that our example in Proposition \ref{ctr_ex}(i)  borrowed an idea from \cite[Example 3.5]{ABGK}. 
In  
\cite{GPZ}, they consider controlled random walks, namely 
random walks that are martingales with uniformly bounded increments and nontrivial jump probabilities
(that may depend on the behavior of the random walks), 
and show that anomalous behavior of the heat kernels can occur in the framework. 
The readers may find further related works in the references of the above papers. 

\subsection{Framework and main results}
Let $\bG=(V,E)$ be a 
connected graph with bounded degree. 
Assume that for each $t\ge 0$, the
graph $\bG$ is endowed with a conductance (weight) $\mu^{(t)}(x,y)$ which is a symmetric nonnegative deterministic function on $V\times V$ such that $\mu^{(t)}(x,y)>0$ if and only if $\{x,y\}\in E$. Suppose 
further that the map $t \mapsto \mu^{(t)}(x,y)$ is right continuous and has left limit
(\abbr{rcll} for short) for each $\{x,y\}\in E$. We call $(\bG, \{\mu^{(t)}(x,y)\})$ 
a time-dependent weighted graph. Let $\mu^{(t)}(x):=\sum_y\mu^{(t)}(x,y)$ for each $x$
and define a measure $\mu^{(t)}$ on $V$ by setting $\mu^{(t)}(A)=\sum_{x\in A}\mu^{(t)}(x)$
for each $A\subset V$. Let $\nu$ be a uniform measure on $V$, that is 
$\nu(A)=|A|$ for $A\subset V$ where $|A|$ is a cardinality of $A$. 
Throughout the paper, we assume the following: there exists  
$c_1\in (0,1]$ such that
\begin{eqnarray}\label{eq:p0-unifcond}
c_1\le \mu^{(t)}(x,y)\le c_1^{-1},\,\qquad\qquad~~\forall \{x,y\}\in E.
\end{eqnarray}

We now define a quadratic form on $(\bG, \{\mu^{(t)}(x,y)\})$ as follows:
\[
\cE_t(f,g)=\frac 12 \sum_{x,y\in V}(f(x)-f(y))(g(x)-g(y))\mu^{(t)}(x,y)
\]
for each $f,g\in H_t^2$, where
\[
H_t^2=\{f:V\to\bR: \sum_{x,y\in V}(f(x)-f(y))^2\mu^{(t)}(x,y)<\infty\}.
\]
Define discrete Laplace operators as follows:
 \begin{eqnarray*} 
  {\cL}_t^C f(x)=
 \sum_y(f(y)-f(x)) \frac {\mu^{(t)}(x,y)}{\mu^{(t)}(x)},~~~
{\cL}_t^V f(x)=\sum_y(f(y)-f(x))\mu^{(t)}(x,y).\end{eqnarray*}
For each $f,g$ that has finite support, we have 
\[\cE_t(f,g)=-(\cL_t^Vf,g)_\nu=-(\cL_t^Cf,g)_{\mu^{(t)}},\]
where 
$(f,g)_\theta:=\sum_{x}f(x)g(x)\theta(x)$ 
for a measure $\theta$.\\  

We next provide definitions for discrete time and continuous time constant/variable speed random walks on $(\bG, \{\mu^{(t)}(x,y)\})$. One way to construct such processes is through the theory  of time-dependent Dirichlet forms (see \cite{O}), but this will require some knowledge of probabilistic potential theory 
and some more notation. 
Here we give a more direct definition.    
For $x,y\in V$, we define 
\[
P^{(t)}(x,y):=\mu^{(t)}(x,y)/\mu^{(t)}(x).
\]

\begin{defn}\label{defn:dtrw}
(i) The $V$-valued stochastic process $\{X_t\}_{t\in\mathbb{N}}$ is called  
a discrete time random walk  on $\bG$, if its transition probabilities at time $t\in\mathbb{N}$
are given by $P(t,x;t+1,y)=P^{(t)}(x,y)$, for any $\{x,y\}\in E$.\\
(ii) The $V$-valued stochastic process $\{Y_t\}_{t\in\mathbb{R}_+}$ of \abbr{RCLL} sample path 
$t \mapsto Y_t$ is called  
a constant speed random walk (in short \abbr{CSRW}), if it waits 
i.i.d. $\exp(1)$ times between successive jumps, and if 
$Y_{T^-} = x$ just prior to the current random jump time $T$, then 
the process jumps across each $\{x,y\} \in E$ with probability 
$P^{(T)}(x,y)$.\\
(iii) The $V$-valued stochastic process $\{Y_t\}_{t\in\mathbb{R}_+}$ of \abbr{RCLL} sample path 
$t \mapsto Y_t$ is called  
a variable speed random walk (in short \abbr{VSRW}), if the holding time of the particle at 
$x\in V$ at time $t\in\mathbb{R}_+$ is independent with the law 
$\exp(\mu^{(t)}(x))$, and if 
$Y_{T^-} = x$ just prior to the current random jump time $T$, then 
the process jumps across each $\{x,y\} \in E$ with probability $P^{(T)}(x,y)$.
\end{defn}

\abbr{CSRW} and \abbr{VSRW} as defined above are associated with the operators $\cL^C_t$ and $\cL_t^V$ respectively, in particular the heat kernel of \abbr{CSRW} $P(Y_t=y|Y_s=x)/\mu^{(t)}(y)$ solves $\frac{\partial}{\partial t}p(s,x;t,y)=\cL^C_tp(s,x;t,y)$, whereas the heat kernel of \abbr{VSRW} $P(Y_t=y|Y_s=x)$ solves $\frac{\partial}{\partial t}p(s,x;t,y)=\cL^V_tp(s,x;t,y)$, where the operators act on the $y$ variable. This can be seen from the fact that independently at time $t$ and vertex $x$, the exponential rate at which the random walk jumps across edge $(x,y)$ is given by $\mu^{(t)}(x,y)/\mu^{(t)}(x)$ for \abbr{CSRW} and $\mu^{(t)}(x,y)$ for \abbr{VSRW}.
\medskip

We first show that, for the \abbr{VSRW} we have the stability of Gaussian heat kernel estimates as 
expected. While we could not find out the precise statement as given below, the proof is a 
careful line by line 
modifications of the known proof (such as the proof in \cite{De}). 
Once again we note that 
it is proved in \cite[Sect. 4]{DD} and \cite[Appendix B]{GOS} that any \abbr{VSRW} on 
${\mathbb Z}^d$ among uniformly elliptic time-dependent conductances 
must satisfy the two-sided Gaussian heat kernel bounds. In the setting of 
time-dependent local regular 
Dirichlet forms on metric measure spaces, similar results are given in 
\cite{St1} and the equivalence of the parabolic Harnack inequalities and 
the volume doubling property plus the Poincar\'e inequalities are given in  \cite{St2}.

 \begin{thm}\label{thm:PIVDPHI}  
Let $\bG=(V,E)$ be 
a connected graph with bounded degree and $\{\mu^{(t)} (x,y): x,y\in V\}$ be 
time-dependent conductances that satisfy \eqref{eq:p0-unifcond}. 
Then the following are equivalent:\\
(a) The graph $\bG$ satisfies 
the volume doubling with constant $C_1<\infty$, namely \begin{align}\label{eq:vvd}
\nu(\bB(x,2r))\le C_1\nu(\bB(x,r))
\end{align}
for all $x\in\bG$, $r>0$; and 
the Poincar{\'{e}} inequality 
holds 
with a constant $C_2<\infty$, namely 
\begin{align}\label{eq:PoinI}
\sum_{x\in \bB(x_0, r)}|f(x)-f_{\bB}|^2\le C_2 r^2\sum_{x,y\in \bB(x_0,2r)}(f(x)-f(y))^2\mu^{(t)}(x,y),
\end{align}
for all $f: V\to \bR$, $x_0\in\bG$, $r>0$, where $f_{\bB}=\sum_{x\in \bB(x_0,r)}f(x)/\nu(\bB(x_0,r))$.\\
(b) The parabolic Harnack inequality holds 
for all non-negative solutions of equation
\begin{align*}
\frac{\partial u(t,x)}{\partial t}={\cL}_t^Vu(t,x).
\end{align*}
That is, set $\eta\in(0,1)$ and $0<\theta_1<\theta_2<\theta_3<\theta_4$, 
we have for all $x_0, s, r$, 
and every non-negative solution on cylinder $Q=[s,s+\theta_4 r^2]\times\bB(x_0,r)$,
\begin{align*}
\sup_{Q_{-}}u\le C_3\inf_{Q_{+}}u,
\end{align*}
where $Q_{-}=[s+\theta_1r^2, s+\theta_2r^2]\times\bB(x_0,\eta r)$ and $Q_{+}=[s+\theta_3r^2, s+\theta_4r^2]\times \bB(x_0,\eta r)$, with some $C_3<\infty$.\\
(b${}^*$) The parabolic Harnack inequality holds for all non-negative solutions of equation
\begin{align*}
\frac{\partial u(t,x)}{\partial t}={\cL}^Vu(t,x), 
\end{align*}
where ${\cL}^V f(x):=\sum_{y: \{x,y\}\in E}(f(y)-f(x))$.\\
(c) The following two-sided heat kernel estimates hold for 
the corresponding \abbr{VSRW}: 
there exist positive constants $C_4,C_5, c_6,c_7<\infty$ such that
\begin{align}\label{eq:GaussUHK}
p(0,x;t,y)
\le \left\{\begin{array}{ll}
\frac{C_4}{\nu(\bB(x,t^{1/2}))}\exp\Big(-C_5 
d(x,y)(1\vee \log(\frac{d(x,y)}{t})\Big), ~~\qquad~\forall t\le d(x,y),\\ 
\frac{C_4}{\nu(\bB(x,t^{1/2}))}\exp\Big(-C_5\frac{d(x,y)^2}{t}\Big), ~~~~~~~~~\qquad\qquad\qquad\qquad~\forall t\ge d(x,y),\end{array}\right.\\
\frac{c_6}{\nu(\bB(x,t^{1/2}))}\exp\Big(-c_7\frac{d(x,y)^2}{t}\Big)\le p(0,x;t,y),~~\qquad~\forall t\ge d(x,y),
\label{eq:GaussLHK}
\end{align}
for all $x,y\in V$, $t>0$ (with a restriction $ t\ge d(x,y)$ in \eqref{eq:GaussLHK}). 
\end{thm}
\begin{remark}
We note that the uniformly elliptic condition \eqref{eq:p0-unifcond} is a natural assumption 
when discussing the stability of \abbr{VSRW}s even for the time-independent case. Indeed, consider 
one parameter family of conductances $\{M\mu (x,y): x,y\in V\}$ with $M\ge 1$ and assume that \eqref{eq:vvd} and \eqref{eq:PoinI} hold when $M=1$. Then they also hold for all $M\ge 1$. However, the corresponding 
\abbr{VSRW} is a constant ($M$ times) time change of the process for $M=1$, so \eqref{eq:GaussUHK} and \eqref{eq:GaussLHK} cannot hold uniformly. The lower bound of \eqref{eq:p0-unifcond} can be deduced by  
applying \eqref{eq:PoinI} with $r=1$ and $f(z)=\delta_{\{x\}}(z)$.
\end{remark}
Since the proof of Theorem \ref{thm:PIVDPHI} is similar to that of the time-independent case, we will simply give a sketch of the proof in the next section. As a consequence of this theorem, we can see that the \abbr{VSRW} on ${\mathbb Z}^d$ among uniformly elliptic time-dependent conductances enjoys the Gaussian heat kernel estimates \eqref{eq:GaussUHK}
and \eqref{eq:GaussLHK}.

\bigskip

In contrast to the above theorem, for the discrete time random walk and the \abbr{CSRW}, one can construct 
a transient random walk on ${\mathbb Z}$ among uniformly elliptic time-dependent conductances 
as in the next 
proposition (cf. \cite[Example 3.5]{ABGK}).

Let $\gamma <1$. We say a time-dependent discrete time random walk is $\gamma$-lazy if $P^{(t)}(x,x)\ge \gamma$ for all $x\in V$ and all $t\ge 0$. 

\begin{ppn}\label{ctr_ex}
(i) For any $\gamma<1$ and $\vep>0$ there exist time-dependent conductances 
$\{\mu^{(t)}(x,x \pm 1), \mu^{(t)}(x,x): x \in \bZ\}$ on $\bZ$ with 
\[
1-\vep\le \mu^{(t)}(x,x \pm 1)\le 1+\vep, \,\qquad~~ \forall x\in \bZ,\, t\in\bN
\]
such that the corresponding discrete time random walk $\{X_t\}_{t\in\bN}$
is $\gamma$-lazy, and it is ballistic and transient almost surely (i.e. it returns to starting point finitely often).

(ii) For any $\vep, c>0$, there exist time-dependent, piecewise-constant conductances 
$\{\mu^{(t)}(x,x \pm 1): x \in \bZ\}$ on $\bZ$ with 
\[
1-\vep\le \mu^{(t)}(x,x \pm 1)\le 1+\vep,  \,\qquad~~ \forall x\in \bZ, \, t\in\bR_+
\]
and the times $\{t_n\}$ at which the conductances change satisfying $t_n/n\to1/c$, 
such that the corresponding \abbr{CSRW} $\{Y_t\}_{t\in\bR_+}$ is ballistic and transient almost surely. 

In particular, both walks violate the Gaussian 
heat kernel on-diagonal lower bound as well as off-diagonal upper bound on $\bZ$. 
\end{ppn}

In this proposition, we give examples where the edge conductances are periodically fluctuating in time. 
It remains open whether adding monotone condition on the edge conductances would recover the expected Gaussian lower bound.

In a similar manner, we also give in the next proposition examples of discrete time random walks and \abbr{CSRW} on $\bZ^2\times\bZ_{\ge0}$ with uniformly elliptic time dependent conductances that violate the on-diagonal Gaussian upper bound. Let $\{e_1,e_2,e_3\}$ be the Cartesian standard basis of $\bZ^3$.

\begin{ppn}\label{ctr_ex_2}
(i) For any $\gamma<1$ and $\ep>0$ there exist time-dependent conductances $\{\mu^{(t)}(x,x\pm e_i),\mu^{(t)}(x,x): x\in\bZ^2\times\bZ_{\ge 0}, i=1,2,3\}$ on $\bZ^2\times\bZ_{\ge0}$ with 
\begin{align*}
1-\ep\le\mu^{(t)}(x,x\pm e_i)\le 1+\ep, \quad  i=1,2,3,\,\forall x\in\bZ^2\times\bZ_{>0}, \, t\in\bN
\end{align*}
such that the corresponding discrete time random walk $\{X_t\}_{t\in\bN}$ is $\gamma$-lazy and recurrent almost surely (i.e. it returns to starting point infinitely often).

(ii) For any $\ep,c>0$, there exist time-dependent, piecewise-constant conductances $\{\mu^{(t)}(x,x\pm e_i): x\in\bZ^2\times\bZ_{\ge 0}, i=1,2,3\}$ on $\bZ^2\times\bZ_{\ge0}$ with 
\begin{align*}
1-\ep\le\mu^{(t)}(x,x\pm e_i)\le 1+\ep, \quad  i=1,2,3, \,\forall x\in\bZ^2\times\bZ_{\ge0},  \, t\in\bR_+
\end{align*}
and the times $\{t_n\}$ at which the conductances change satisfying $t_n/n\to1/c$, such that the corresponding \abbr{CSRW} $\{Y_t\}_{t\in\bR_+}$ is recurrent almost surely.

In particular, both walks violate the Gaussian heat kernel on-diagonal upper bound on $\bZ^2\times\bZ_{\ge0}$.
\end{ppn}

As a consequence, for the random walks in Propositions \ref{ctr_ex} and \ref{ctr_ex_2} the stability of Gaussian heat kernel estimates of Theorem \ref{thm:PIVDPHI} does not hold, even though (VD) and (PI(2)) hold uniformly.

\end{section}

\begin{section}{Proof.} 
{\emph{Sketch of the proof of Theorem \ref{thm:PIVDPHI}.}}

(a) $\Rightarrow$ (b) $\Rightarrow$ (c): 
As explained in \cite[pg 374-375]{DD}, it is possible to adapt Delmotte's argument (\cite{De}) here by setting, in their notation, $\mu_{xy}=a(t,x,y)=: \mu^{(t)}(x,y)$ and $m(x)=:1$. 
(Note that Delmotte's proof is for the discrete time random walk and \abbr{CSRW}.) 
Note that although the framework of \cite{DD} is $\bZ^d$, the same modification can be employed for general $\bG$. In both \cite{De,DD} the term 
${\cE}(t,D)=\exp(-D\arg\sinh\frac{D}{t}+t(\sqrt{1+D^2/t^2}-1))$ appears in the off-diagonal bounds, 
but simple computations show that it is comparable with the exponential parts of 
\eqref{eq:GaussUHK}, \eqref{eq:GaussLHK}. 
Let us now overview the proof. Assuming (a), \eqref{eq:p0-unifcond} translates to a Poincar{\'{e}} inequality that holds uniformly for all $t$, as well as a weighted Poincar{\'{e}} inequality and a Sobolev-Poincar{\'{e}} inequality needed along the way (\cite[Proposition 2.2, 2.4]{De}). Since \cite[Section 2]{De} is itself in continuous time, one thus can re-produce the entire section, resulting in a parabolic Harnack inequality (b). Now (b) implies the on-diagonal upper bound and the near diagonal (i.e. for $d(x,y)^2\le t$) lower bound for the heat kernel of $\{X_t\}$ and its dual process. (Note that unlike the time-independent case, $p(0,x,t,y)$ is 
no longer equal to $p(0,y,t,x)$. However, it holds that $p(0,x,t,y)=p^*(0,y,t,x)$ where 
$p^*(\cdot,\cdot,\cdot,\cdot)$ is the heat kernel for the time reversal conductances, i.e. 
$\{\mu^{(t-\cdot)}(x,y)\}$; cf. \cite[Lemma 1.5]{St1}.)  
The off-diagonal upper bound can be deduced from the on-diagonal one and the integrated maximum principle 
using the Davies' argument. The off-diagonal lower bound (in the range $d(x,y)\le t$) follows from the 
near diagonal one by the usual chain argument. See \cite[Section 3.1]{De} for details. 

(a) $\Leftrightarrow$ (b${^*}$): Note that \eqref{eq:PoinI} is equivalent to the inequality 
where the right hand side is changed to 
$C_2' r^2\sum_{x,y\in \bB(x_0,2r)}(f(x)-f(y))^2$. 
So, the equivalence for time-independent case (that can be proved similarly to \cite{De}) implies the desired equivalence. 

(b) $\Rightarrow$ (b${^*}$): This is trivial.

(c) $\Rightarrow$ (b): 
This can be proved using the Balayage argument as in 
\cite[Theorem 3.10]{De} (see the proof of \cite[Theorem 1.5]{BKM} for more details on 
the Balayage argument in the setting of continuous time Markov chains).
In order to apply the Balayage argument, the existence of the space-time dual process is required --in this case, we know the existence by 
using the time reversal conductances mentioned above. In the proof, we need the following estimate
\[
\sup_{0<s\le R^2}p(0,x,s,y)\le \frac{c_1}{\nu(\bB(x_0,R))}
\,~~~\qquad\mbox{for all }~x,y\in \bB(x_0,2R)\,~~\mbox{with }~
d(x,y)\ge R,\]
which can be deduced by \eqref{eq:GaussUHK} and the fact that
there exists $\beta>0$ such that 
$\nu(\bB(x,R))\le cR^\beta$ for all $x\in V, R\ge 1$.
The last inequality is a consequence of \eqref{eq:vvd} and 
the degree of the graph being bounded. 
We note that we do not need the heat kernel lower bound for $t\le d(x,y)$ to establish (b). 
\qed

\bigskip 
{\emph{Proof of Proposition. \ref{ctr_ex}.}}
(i) Here $\mathbb{G}=\bZ$ and we set the edge conductances to be 
\begin{align*}
\mu^{(t)}(i,i-1)=1-\vep, \,\mu^{(t)}(i,i)=b,\, \mu^{(t)}(i,i+1)=1+\vep, \,\text{when } t+i\text{ is even};\\
\mu^{(t)}(i,i-1)=1+\vep,\, \mu^{(t)}(i,i)=b', \,\mu^{(t)}(i,i+1)=1-\vep, \,\text{when }
t+i \text{ is odd}\,,
\end{align*}
with $b/(b+2)=\gamma$ and $b'/(b'+2)=\gamma'>\gamma$. We start at
$X_0=0$ and notice that this random walk has two possible states:
either $X_t$ is at state $A_+$ with his right edge having 
conductance $1+\vep$, or it is at state $A_-$ with his right
edge having conductance $1-\vep$. Whenever the random walk $X_t$ 
moves either to its left or right vertex, it keeps the current state, while 
if it stays put (i.e. $X_{t+1}=X_t$), then due to the change of conductance 
values, it moves to the opposite state. Let $\{Z_t\}_{t\in\bN}$ be the $\{A_{\pm}\}$-valued Markov chain describing the state of $\{X_t\}$, then
the transition probabilities  of $Z_t$ are thus 
\begin{align}
q(A_+,A_+)=1-\gamma,\;\; q(A_+,A_-)=\gamma,\;\; q(A_-,A_-)=1-\gamma',\;\; 
q(A_-,A_+)=\gamma' \,.  \label{matrix}
\end{align}
and its invariant measure is
\begin{align}
\pi(A_+)=\frac{\gamma'}{\gamma'+\gamma},\,\quad 
\pi(A_-)=\frac{\gamma}{\gamma'+\gamma} \,,  \label{inv-meas}
\end{align}
whereas by the strong law for occupation time $N_t(\cdot):=\sum_{i=0}^{t-1}\bI_{\{Z_i=\cdot\}}$ (Cf. \cite[(5.5) pg 320]{Du}),
\begin{align}
N_t(A_{\pm})/t\stackrel{a.s.}{\to}\pi(A_{\pm}). \label{occup-time}
\end{align}
Further, whenever at state $A_+$ the random walk has drift 
$\Delta(A_+)=\vep(1-\gamma)$ to its right while
at state $A_-$ it has drift $\Delta(A_-)=- \vep(1-\gamma')$. We enumerate sequentially the random times $m_1<m_2< ...$ when the random walk is at state $A_+$, and similarly enumerate the random times $n_1<n_2<...$ when the random walk is at state $A_-$, then $\cS_+:=\{D_i:=X_{i+1}-X_i,\, i\in\{m_1,m_2,...\}\}$ are i.i.d. with drift $\Delta(A_+)$, and $\cS_-:=\{D_i:\,i\in\{n_1,n_2,...\}\}$ are i.i.d. with drift $\Delta(A_-)$, while $\cS_{\pm}$ are also mutually independent. Hence by the strong law of large numbers (\abbr{SLLN}) and (\ref{occup-time}), we have that 
\begin{align}
\frac{X_t}{t}&=\frac{\sum_{i=0}^{t-1}D_i\bI_{\{D_i\in\cS_+\}}}{N_t(A_+)}\frac{N_t(A_+)}{t}+\frac{\sum_{i=0}^{t-1}D_i\bI_{\{D_i\in\cS_-\}}}{N_t(A_-)}\frac{N_t(A_-)}{t}\nonumber\\
&\stackrel{a.s.}{\to}
\Delta(A_+)\pi(A_+)+\Delta(A_-) \pi(A_-)
=\vep \frac{\gamma'(1-\gamma)-\gamma (1-\gamma')}{\gamma'+\gamma} 
= \vep \frac{\gamma'-\gamma}{\gamma'+\gamma}=:\beta > 0.
\label{eq:balli}\end{align}
It is thus ballistic and transient almost surely, violating the Gaussian heat kernel on-diagonal lower bound (\ref{eq:GaussLHK}).
If the Gaussian 
off-diagonal upper bound (\ref{eq:GaussUHK}) holds, then 
integrating over the region $y\in[(\beta-\ep)t,(\beta+\ep)t]$ for any $\ep\in(0,\beta)$, we see that 
$\bP(X_t\in [(\beta-\ep)t,(\beta+\ep)t])$ decays 
exponentially for $t$, which contradicts \eqref{eq:balli}. So this walk 
also
violates the Gaussian heat kernel off-diagonal upper bound (\ref{eq:GaussUHK}).  

To have a non-lazy example, set $b=b'=0$ and observe that $\{X_t\}$ then keeps the state $A_+$ at all times.\\

(ii) Here again $\bG=\bZ$, and let $\{\tau_k\}_{k\in\bN}$ be the successive jump times of a Poisson process of intensity $c-1\in(0,\infty)$, with $\tau_0=0$, independent of the \abbr{CSRW} $\{Y_t\}$, and then we set the edge conductances to be 
\begin{align*}
\mu^{(t)}(i,i+1)=1-\ep, \quad \mu^{(t)}(i+1,i+2)&=1, \quad \mu^{(t)}(i+2,i+3)=1+\ep, \\
 &\text{when }t\in[\tau_k,\tau_{k+1}), \text{ and }i\equiv k \text{ mod } 3.
\end{align*}
We start at $Y_0=0$ and notice that this \abbr{CSRW} has three possible states: either it is at state $A_1$ with left/right (\abbr{L/R} for short) edge conductances $1+\ep$, $1-\ep$, or it is at state $A_2$ with \abbr{L/R} edge conductances $1-\ep$, $1$, or at state $A_3$ with \abbr{L/R} edge conductances $1$, $1+\ep$. On the other hand, there are two independent Poisson clocks, one ($``\sC_E"$) governing the environment shift which has intensity $c-1$, and one ($``\sC_J"$) governing jumps of \abbr{CSRW} which has intensity $1$. Denote $\{T_k\}_{k\in\bN}$ the sequence of times when 
the state of $\{Y_t\}$ changes, then it is the successive jump times of a Poisson process of intensity $c\in(1,\infty)$. Let $\{Z_k\}_{k\in\bN}$ be the $\{A_1,A_2,A_3\}$-valued process describing the state of $Y_{T_k}$, then the transition from $Z_k$ to $Z_{k+1}$ is determined by which clock rings first, and in case $\sC_J$ does, what are the adjacent edge conductances, but not on $\{Z_0,...,Z_{k-1}\}$. In other words, the process $\{Z_k\}$ is a time-homogeneous Markov chain with state space $\{A_1,A_2,A_3\}$.

Using properties of exponential distribution (i.e. if $\xi_1$ and $\xi_2$ are independent exp($\gamma_1$) and exp($\gamma_2$) random variables, 
then 
$\bP(\xi_1<\xi_2)=\gamma_1/(\gamma_1+\gamma_2)$), one can calculate the transition probabilities of $\{Z_k\}$:
\begin{align*}
q(A_1,A_2)=\frac{1-\ep}{2c}, \, q(A_1,A_3)=&1-\frac{1-\ep}{2c},\, q(A_2, A_3)=\frac{1}{(2-\ep)c},\,
q(A_2, A_1)=1-\frac{1}{(2-\ep)c},\\
q(A_3,A_1&)=\frac{1+\ep}{(2+\ep)c},\,  q(A_3,A_2)=1-\frac{1+\ep}{(2+\ep)c}.
\end{align*} 
and its invariant measure is proportional to 
\begin{align*}
\pi=
\big[\, 2[(-4c^2+2c-1)+&(c-1)\ep+c^2\ep^2],\, (2-\ep)[(4c^2-2c+1)+(2c^2-2c)\ep-\ep^2],\\
& (2+\ep)[(4c^2-2c+1)+(-2c^2+3c-1)\ep-c\ep^2]\, 
\big].
\end{align*}
Further, the drift $Y_t$ is subject to when at states $A_i$, $i=1,2,3$ for its immediate next change of state is 
$\Delta=\big[-\frac{\ep}{c}, \frac{\ep}{(2-\ep)c}, \frac{\ep}{(2+\ep)c}\big]$.
By \abbr{SLLN} the speed of $\{Y_t\}$ is proportional to $\pi\cdot\Delta$, and one can check that the speed is positive when $\ep\in(-1,-\frac{3}{2c+1})\cup(0,1)$, and negative when $\ep\in(-\frac{3}{2c+1},0)$ (see also Remark \ref{speed}). This implies that for arbitrary $\ep\in(0,1)$, $c>1$, $\{Y_t\}$ has non-zero speed, w.p.1 under the annealed measure on the environment.
 
Furthermore, the annealed result implies that for a.e. realization of the isolated Poisson jump times $\{\tau_n\}$, $\{Y_t\}$ is w.p.1. 
transient, ballistic and violates the Gaussian heat kernel off-diagonal upper bound (\ref{eq:GaussUHK})
in the quenched sense. That is, there exist some (in fact uncountably many) choices of non-random $\{t_n\}$ with $t_n/n\to1/(c-1)$, so that with the conductances changing at times $\{t_n\}$ the corresponding \abbr{CSRW} on $\bZ$ is transient, ballistic and violates the off-diagonal upper bound (\ref{eq:GaussUHK}).
\qed
\begin{remark}\label{speed}
An intuitive explanation of the phenomenon concerning the region of positive/negative speed in the above example is as follows. Its asymmetry comes from the fact that, although the conductances are symmetric in both directions,  the environment is shifting only to the right, and this breaks the symmetry of $\ep$. Also, when $\ep$ is sufficiently close to $-1$, the speed becomes positive again. Take the special case when $c\to\infty$, then the shift of conductances is so quick that at every time the \abbr{CSRW} jumps (which happens independently at rate $1$), its neighborhood is one of the three choices with almost equal probability. But the drift at these three neighborhoods are $-\ep$, $\ep/(2-\ep)$, $\ep/(2+\ep)$ respectively, one can check that their average is positive regardless of $\ep\in(-1,0)\cup(0,1)$.
\end{remark}
\begin{remark}\label{osc}
The effect of oscillating edge conductances can be mapped to monotone but unboundedly increasing or decreasing conductances. Take one dimension and discrete time for example, with $a>0$ the oscillating conductances $\{..,1,a,1,a...\}$ on $\bZ$ shifting at speed $1$, is equivalent to setting at time $2n$ the conductances $\{.., a^{2n},a^{2n+1},a^{2n},a^{2n+1}...\}$ and at time $2n+1$ $\{.., a^{2n+2},a^{2n+1},a^{2n+2}, a^{2n+1}..\}$ etc. However, we expect that among monotone and uniformly elliptic conductances, the walks follow recurrence/transience of the starting and ending graphs. This has been proved for discrete time non-lazy walks on trees in \cite[Theorems 5.1, 5.2]{ABGK}. 
\end{remark}

{\emph{Proof of Proposition \ref{ctr_ex_2}.}} 
(i) Here the vertex set is $\V=\bZ^2\times\bZ_{\ge0}$ and we set the edge conductances to be 
\begin{align*}
&\mu^{(t)}(\xi,\xi+e_3)=1+\ep, \, \mu^{(t)}(\xi,\xi-e_3)=1-\ep,\, \mu^{(t)}(\xi,\xi)=b, \\
&\quad\quad\text{ when }\xi=(i,j,k)\in\bZ^2\times\bZ_{>0}, \, t+i+j+k \text{ is odd};\\
&\mu^{(t)}(\xi,\xi+e_3)=1-\ep, \, \mu^{(t)}(\xi,\xi-e_3)=1+\ep,\, \mu^{(t)}(\xi,\xi)=b', \\
&\quad\quad\text{ when }\xi=(i,j,k)\in\bZ^2\times\bZ_{>0}, \, t+i+j+k \text{ is even};\\
&\mu^{(t)}(\xi,\xi\pm e_l)=1, \, l=1,2,\text{ for all }\xi\in\bZ^2\times\bZ_{>0}, \text{ and all }t\, ;\\
&\mu^{(t)}(\xi,\xi\pm e_l)=0, \, l=1,2, \text{ for all }\xi=(i,j,0), \text{ and all }t, \,\\
&\mu^{(t)}(\xi,\xi)=f,\, \text{when }\xi=(i,j,0),\, t+i+j\text{ is odd};\\
&\mu^{(t)}(\xi,\xi)=f',\,  \text{when }\xi=(i,j,0),\, t+i+j\text{ is even.}
\end{align*}
with $b/(b+6)=f/(f+1+\ep)=\gamma$ and $b'/(b'+6)=f'/(f'+1-\ep)=:\gamma'<\gamma$.

Starting at $X_0=\underline{0}$ the random walk has two possible states. Either it is at state $A_+$ with upper edge conductance $1+\ep$, or it is at state $A_-$ with upper edge conductance $1-\ep$. 
Whenever the random walk moves, it keeps the same state; and whenever it stays put, it changes to the opposite state. Let $\{Z_t\}_{t\in\bN}$ denote the state of $\{X_t\}$. Define the sequence of stopping times $\{\sigma_n\}$ starting from $\sigma_0=0$, and for $i\ge 1$, $\sigma_i:=\inf\{t>\sigma_{i-1}: (X_t)_3=0\}$, and let $M_n:=\big((X_{\sigma_n})_1,(X_{\sigma_n})_2\big)$ be the two-dimensional random walk on $\bZ^2\times\{0\}$. When $R_t:=(X_t)_3>0$, the state transition probabilities $\{q(\cdot,\cdot)\}$ are given by (\ref{matrix}) and they have an invariant measure $\pi(\cdot)$ given by (\ref{inv-meas}); whereas at state $A_+$, the random walk has drift $\Delta(A_+)=(2\ep)/(6+b)$ and at state $A_-$ it has drift $\Delta(A_-)=-(2\ep)/(6+b')$.  Let 
\begin{align*}
\beta:=\Delta(A_+)\pi(A_+)+ \Delta(A_-)\pi(A_-)=\frac{2\ep}{6+b}\frac{\gamma'}{\gamma+\gamma'}-\frac{2\ep}{6+b'}\frac{\gamma}{\gamma+\gamma'}<0.
\end{align*}
Because $\beta<0$, by the large deviation arguments, 
there exists some positive constant $c_1=c_1(\beta)$ such that for all $k$ large enough and every $n$,
\begin{align}
\bP(||D_n||>\sqrt{2}k|\cF_n)\le\bP(\sigma_{n+1}-\sigma_n>k|\cF_n)=\bP(\min_{1\le i\le k}R_{\sigma_n+i}>0|\cF_n)\le c_1^{-1}e^{-c_1k}, \label{tail}
\end{align}
where $\cF_n:=\cF^X_{\sigma_n}$ is the canonical filtration of $\{X_t\}$ stopped at $\sigma_n$, and $D_n:=M_{n+1}-M_n$. We enumerate sequentially the random times $m_1<m_2<...$ when the state of $M_n$ is $A_+$, and similarly the random times $n_1<n_2<...$ when the state of $M_n$ is $A_-$. Then the collection $\cS_+:=\{D_i:\, i\in\{m_1,m_2...\}\}$ are i.i.d. with some law $\nu_+$, and the collection $\cS_-:=\{D_i:\, i\in\{n_1,n_2...\}\}$ are i.i.d. with some law $\nu_-$, while $\cS_{\pm}$ are mutually independent. Also, the sequence of states $\{Z_{\sigma_n}\}$ approach an invariant measure which we denote by $\widetilde{\pi}(\cdot)$ (different from $\pi(\cdot)$). Let $p(\cdot,\cdot)$ be the heat kernel of $\{M_n\}$, and $p^{\pm}(\cdot,\cdot)$ the heat kernels of aperiodic random walks of i.i.d. increments with law $\nu_{\pm}$, which have all moments finite by (\ref{tail}), 
we have 
by the strong law for state occupation time of $\widetilde{\pi}(\cdot)$ that 
\begin{align}
\lim_{n\to\infty}\frac{p(n,\underline{0})}{\big(p^+(\widetilde{\pi}(A_+)n,\cdot)*p^-(\widetilde{\pi}(A_-)n,\cdot)\big)(\underline{0})}=1, \label{hk-est}
\end{align}
where $*$ denotes convolution on $\bZ^2$. By the local central limit theorem for $p^{\pm}(\cdot,\cdot)$ (\cite[Theorem 2.3.5]{LL}), there exists some positive constant $c_2$ such that for all $z$ satisfying $||z||\le\sqrt{n}$, we have that $p^{\pm}(\widetilde{\pi}(A_{\pm})n,z)\ge c_2/n$, therefore
\begin{align}
\big(p^+(\widetilde{\pi}(A_+)n,\cdot)*p^-(\widetilde{\pi}(A_-)n,\cdot)\big)(\underline{0})\ge\int_{\{||z||\le\sqrt{n}\}}p^+(\widetilde{\pi}(A_+)n,z)p^-(\widetilde{\pi}(A_-)n,-z)dz\ge c_2^2/n, \nonumber\\
\label{local-clt}
\end{align}
which by (\ref{hk-est}) then implies that $\{M_n\}$ is recurrent almost surely, and the same for $\{X_t\}$.  As a consequence, the Gaussian heat kernel on-diagonal upper bound (\ref{eq:GaussUHK}) does not hold for $\{X_t\}$. 

To have a non-lazy example, set $b=b'=f=f'=0$ and observe that $\{X_t\}$ then keeps the state $A_-$ at all times.\\

(ii) Here again $\V=\bZ^2\times\bZ_{\ge0}$, and let $\{\tau_n\}_{n\in\bN}$ be the successive jump times of a Poisson process of intensity $c-1\in(0,\infty)$, with $\tau_0=0$, independent of the \abbr{CSRW} $\{Y_t\}$, and then we set the edge conductances to be
\begin{align*}
\mu^{(t)}(\xi,\xi+e_3)=1+\ep,\, &\mu^{(t)}(\xi,\xi-e_3)=1-\ep, \, \mu^{(t)}(\xi,\xi\pm e_l)=1+\ep/2,\, l=1,2\\
&\text{when }\xi=(i,j,k)\in\bZ^2\times\bZ_{>0},\, t\in[\tau_n,\tau_{n+1}), \, n+k \text{ is odd};\\
\mu^{(t)}(\xi,\xi+e_3)=1-\ep,\, &\mu^{(t)}(\xi,\xi-e_3)=1+\ep, \, \mu^{(t)}(\xi,\xi\pm e_l)=1-\ep/2,\, l=1,2\\
&\text{when }\xi=(i,j,k)\in\bZ^2\times\bZ_{>0},\, t\in[\tau_n,\tau_{n+1}), \, n+k \text{ is even};\\ 
\mu^{(t)}(\xi,\xi\pm e_l)=1+\ep/2, \, &\text{when }\xi=(i,j,0), \,l=1,2,\, t\in[\tau_n,\tau_{n+1}), \, n\text{ is odd};\\
\mu^{(t)}(\xi,\xi\pm e_l)=1-\ep/2, \, &\text{when }\xi=(i,j,0), \, l=1,2, \, t\in[\tau_n,\tau_{n+1}),\, n\text{ is even}.
\end{align*} 

Starting at $Y_0=\underline{0}$ the \abbr{CSRW} has two possible states: either it is at state $A_+$ with its upper edge conductance $1+\ep$, or at state $A_-$ with its upper edge conductance $1-\ep$. 
Let $\{T_n\}_{n\in\bN}$ be the sequence of times when the state of $\{Y_t\}$ changes, then it is the successive jump times of a Poisson process of intensity $c\in(1,\infty)$. Let $\{Z_n\}_{n\in\bN}$ be the $\{A_{\pm}\}$-valued time-homogeneous Markov chain describing the state of $\{Y_{T_n}\}$. When $R_{T_n}:=(Y_{T_n})_3>0$, the transition probabilities of $\{Z_n\}$ are given by
\begin{align*}
q(A_+,A_-)&=\frac{c-1}{c}+\frac{1}{(3+\ep)c}, \, q(A_+,A_+)=\frac{2+\ep}{(3+\ep)c},\\
q(A_-,A_+)&=\frac{c-1}{c}+\frac{1}{(3-\ep)c}, \, q(A_-,A_-)=\frac{(2-\ep)}{(3-\ep)c},
\end{align*}
and they have an invariant measure (where $\cZ$ is the normalizing constant)
\begin{align*}
\pi(A_+)=\big(\frac{c-1}{c}+\frac{1}{(3-\ep)c}\big)/\cZ,\quad \pi(A_-)=\big(\frac{c-1}{c}+\frac{1}{(3+\ep)c}\big)/\cZ;
\end{align*}
whereas at state $A_+$, the \abbr{CSRW} has drift (for its immediate next change of state) $\Delta(A_+)=(2\ep)/((6+2\ep)c)$, while at state $A_-$, it has drift $\Delta(A_-)=-(2\ep)/((6-2\ep)c)$.

Hence working analogously to part (i) and defining $\{\sigma_n\}$ for the embedded Markov chain $\{Y_{T_n}\}$,  with
\begin{align*}
\wh{\beta}&:=\Delta(A_+)\pi(A_+)+\Delta(A_-)\pi(A_-)\\
&=\frac{\ep}{(3+\ep)c}\frac{c-1+(3-\ep)^{-1}}{c\cZ}-\frac{\ep}{(3-\ep)c}\frac{c-1+(3+\ep)^{-1}}{c\cZ}<0
\end{align*} 
we thus have the tail bound (\ref{tail}) holding and can carry out the rest of the proof resulting in the heat kernel estimates (\ref{hk-est})-(\ref{local-clt}) as well as almost sure recurrence of $\{Y_t\}$ under the annealed measure on the environment $\{\tau_n\}$. By the same argument as in the proof of Proposition \ref{ctr_ex} (ii), this implies that there exists some non-random $\{t_n\}$ with $t_n/n\to1/(c-1)$, such that almost surely under the quenched measure with the conductances changing at $\{t_n\}$, $Y_t$ is recurrent. It thus follows that the Gaussian heat kernel on-diagonal upper bound does not hold for $\{Y_t\}$.
\qed

\end{section}

\bigskip
\noindent {\bf Acknowledgment. } 
This work was initiated while the second author was visiting 
Stanford University. 
The authors are very grateful to A. Dembo for fruitful discussions and very helpful comments.  
They also thank anonymous referees for the thoughtful comments to the first version of the paper, 
in particular for pointing out an error in Theorem \ref{thm:PIVDPHI}.

\end{document}